\documentclass[12pt]{amsart}

\usepackage{amsfonts,amsthm,amscd,amsmath,latexsym,amssymb}

  \newtheorem{theorem}{Theorem}[section]
  \newtheorem{corollary}[theorem]{Corollary}
  \newtheorem{lemma}[theorem]{Lemma}
  \newtheorem{proposition}[theorem]{Proposition}
  \newtheorem{remark}[theorem]{Remark}
  \newtheorem{mtheorem}[theorem]{Main Theorem}
  \newtheorem{observation}[theorem]{Observation}
  \newtheorem{defn}[theorem]{Definition}
  \theoremstyle{remark}

  \newcommand{\nc}{\newcommand}
  \newcommand{\be}{\begin{equation}}
  \newcommand{\ee}{\end{equation}}
  \newcommand{\bea}{\begin{eqnarray}}
  \newcommand{\eea}{\end{eqnarray}}

\newcommand{\N}{{\mathbb{N}}}
\newcommand{\F}{{\mathbb{F}}}
\newcommand{\lra}{\longrightarrow}

  \newcommand{\bc}{\begin{center}}
  \newcommand{\ec}{\end{center}}
  
  \nc{\bmth}{\begin{mtheorem}} \nc{\emth}{\end{mtheorem}}
  \nc{\bth}{\begin{theorem}}
  \nc{\bpr}{\begin{proposition}} \nc{\epr}{\end{proposition}}
  \nc{\ble}{\begin{lemma}}
   \nc{\ele}{\end{lemma}}
  \nc{\bco}{\begin{corollary}} \nc{\eco}{\end{corollary}}
  \nc{\bre}{\begin{remark}} \nc{\ere}{\end{remark}}
  \nc{\bob}{\begin{observation}} \nc{\eob}{\end{observation}}

  \nc{\f}{\frac} \nc{\rw}{\rightarrow} \nc{\To}{\longrightarrow}
  \nc{\Rw}{\Rightarrow}

  \nc{\nt}{\stackrel{\sim}{\nabla}} \nc{\De}{\Delta}

  \nc{\na}{\nabla} \nc{\al}{\alpha}
   \nc{\bet}{\beta}
   \nc{\va}{\vartheta}
  \nc{\ga}{\gamma} \nc{\G}{\Gamma}
   \nc{\la}{\lambda}  \nc{\ka}{\kappa}
   \nc{\La}{\Lambda}
    \nc{\cS}{\mathcal{S}}
  \nc{\si}{\sigma} \nc{\de}{\delta} \nc{\ep}{\varepsilon}
  \nc{\ei}{\ep _i} \nc{\ti}{\tilde}

    \nc{\ea}{\ep _{\al}}
   \nc{\ej}{\ep _j}
   \nc{\eb}{\ep _{\bet}}
  \nc{\Z}{\mathbb{Z}} \nc{\Te}{\Theta}
  \nc{\om}{\omega}\nc{\Om}{\Omega} \nc{\ro}{\rho}
  \nc{\cG}{\mathcal{G}}\nc{\cR}{\mathcal{R}}
   \nc{\Or}{\mathcal{O}} \nc{\cA}{\mathcal{A}} \nc{\cD}{\mathcal{D}}
  \nc{\C}{{\mathbb{C}}}
  \nc{\q}{\mathbf{Q}}
  \nc{\Q}{\mathbb{Q}}\nc{\bD}{\mathbb{D}} \nc{\bA}{\mathbb{A}}
   \nc{\bP}{\mathbb{P}}
  \nc{\Co}{\mathbb{C}} \nc{\cH}{\mathcal{H}}
  \nc{\cF}{\mathcal{F}}\nc{\cN}{\mathcal{N}}
  \nc{\ca}{\emph{\textbf{a}}} \nc{\bU}{\emph{\textbf{U}}}
\nc{\cM}{\mathcal{M}}  \nc{\cf}{\emph{\textbf{f}}}

  \nc{\bbf}{\mathbf{f}}  \nc{\cX}{\mathcal{X}}
\nc{\Gn}{\mathrm{GL}_n} \nc{\GL}{\mathrm{GL}_2}
\nc{\Sn}{\mathrm{SL}_n} \nc{\SL}{\mathrm{SL}_2}
\nc{\Mn}{\mathrm{M}_n} \nc{\M}{\mathrm{M}_2}

 \newcommand{\cO}{{\mathcal O}}
 \newcommand{\Hom}{{\hbox{\rm Hom}}}

 \newcommand{\symb}{{\hbox{\rm Symb}}}
\newcommand{\bound}{{\hbox{\rm Bound}}}
 \newcommand{\Manin}{{\hbox{\rm Manin}}}
 \newcommand{\Xn}{{\bigl((\Z/p^n\Z)^2\bigr)^\prime}}
 \newcommand \wK {\widetilde K}
 \newcommand \wKn {\wK_2(R_n)}
 \newcommand{\wphi}{{\widetilde\phi}}
 \newcommand{\wpsi}{{\widetilde\psi}}
 \newcommand{\bx}{{\bf x}}
 \newcommand{\bm}{{\bf e}}
 \newcommand{\symbn}{{\symb_{\Gamma_0(p^n)} }}
\newcommand{\boundn}{{\bound_{\Gamma_0(p^n)} }}
\newcommand{\Maninn}{{\Manin_{\Gamma_0(p^n)} }}
\newcommand{\Manininf}{{\Manin^\infty_{\Gamma_0(p^n)} }}
\newcommand{\Gamman}{{\Gamma_0(p^n)}}
\newcommand{\Sigman}{{\Sigma_0(p^n)}}
\newcommand{\pmat}{\left(\begin{matrix}}
\newcommand{\epmat}{\end{matrix}\right)}
\newcommand{\psmat}{\left(\begin{smallmatrix}}
\newcommand{\epsmat}{\end{smallmatrix}\right)}

\begin{document}
\hoffset=-.7in \voffset=-0.05in

\small

  \title[The Steinberg Symbol and Special Values of L-functions]{The Steinberg Symbol and Special Values of L-functions}
  \author{Cecilia Busuioc}

  \address{}

  \date{\today}

  \subjclass{}

 \begin{abstract}

 The main results of this article concern the definition of a
compactly supported cohomology class for the congruence group
$\Gamma_0(p^n)$ with values in the second Milnor $K$-group (modulo
$2$-torsion) of the ring of $p$-integers of the cyclotomic
extension $\Q(\mu_{p^n})$. We endow this cohomology group with a
natural action of the standard Hecke operators and discuss the
existence of special Hecke eigenclasses in its parabolic
cohomology. Moreover, for $n=1$, assuming the non-degeneracy of a
certain pairing on $p$-units induced by the Steinberg symbol when
$(p,k)$ is an irregular pair, i.e. $p|\frac{B_k}{k}$, we show that
the values of the above pairing are congruent mod $p$ to the
$L$-values of a weight $k$, level $1$ cusp form which satisfies
Eisenstein-type congruences mod $p$, a result that was predicted
by a conjecture of R. Sharifi.

\end{abstract}

  \maketitle

 \section{Introduction}

 Let $p^n>1$ be a power of a positive prime $p$, $R_n := \Z[\mu_{p^n},\frac{1}
{p}]$, and $G_n := Gal\left(\Q(\mu_{p^n})/\Q\right)$.   Then $G_n$
acts naturally on Milnor's $K$-group, $K_2^M(R_n)$.  We let the
congruence group $\Gamma_0 := \Gamma_0(p^n)$ act on $K_2^M(R_n)$
via the homomorphism $\Gamma_0 \rightarrow G_n$ sending $\gamma\in
\Gamma_0$ to $\sigma_a$ where $a$ is the upper left hand entry of
$\gamma$ and $\sigma_a\in G_n$ is given by $\sigma_a(\zeta) =
\zeta^a$ for every $\zeta \in \mu_{p^n}$.   In this paper, we will
define a modular symbol
$$
\phi_n \in H^1_c\left(\Gamma_0, \wK_2(R_n)\right)
$$
where $\wK_2(R_n) = K_2^M(R_n)$ (mod $2$-torsion).   We note that
the action of $G_n$ on $R_n$ induces a natural action of $G_n$ on
the above cohomology group. Moreover, we will endow this
cohomology group with a natural action of the Hecke operators
$T_\ell$, $\ell\ne p$, and will prove the following theorem.

 \bth Let $\varphi_n \in H^1_{par}(\Gamma_0,\wK_2(R_n))$ be the image of $\phi_n$ under the
 canonical map.  Then
\begin{enumerate}
\item $\varphi_n |T_2 = (\sigma_2 + 2)\varphi_n$ if $p\ne 2$; and
\item $\varphi_n |T_3 = (\sigma_3 + 3)\varphi_n$ if $p\ne 3$.
\end{enumerate}
\end{theorem}

Now set $n=1$,  let $k \geq 2$ be an even integer and suppose $p >
3$.  Let $R = R_1$, $G=G_1$.  Let $E = R^\times/R^{\times p}$ be
the group of units modulo $p^{th}$ powers of units in $R$. Then we
may decompose $E$ as a direct sum
$$
E = \bigoplus_{i=0}^{p-2} E^{(1-i)}
$$
where $E^{(1-i)}$ denotes the $\F_p$-submodule of $E$ on which $G$
acts via $\omega^{1-i}$ and $\omega : G \lra \F_p^\times$ is the
canonical isomorphism.    It is well-known that $E^{(1)} = \mu_p$,
and $E^{(1-i)} = 0$ for even $i$ satisfying $0< i \le p-3$.  For
$i$ odd, let $\eta_i\in E^{(1-i)} $ be the image of the cyclotomic
$p$-unit $1-\zeta_p$ under the canonical projection $R^\times \lra
E^{(1-i)}$. For $i=1,3,\ldots,p-2$, let $\xi_i
=\{\eta_{k-i},\eta_{i}\}\in(K_2^M(R)/pK_2^M(R))^{(2-k)}.$ We
remark that Vandiver's conjecture implies that the cyclotomic
$p$-units generate $E$, therefore also that the symbols $\xi_i$,
($i=1,\ldots,p-2$) generate $(K_2^M(R)/pK_2^M(R))^{(2-k)}.$  In
what follows, we will often assume the following hypothesis:
\medskip

\noindent {\bf Hypothesis $(H_k)$}:  There exists a non-zero
$G$-equivariant map
$$\rho:  K_2^M(R)/pK_2^M(R) \rightarrow
\F_p(\omega^{2-k}).
$$
and an odd integer $i$ with $1<i<k-1$ such that the pairing
induced by the composition of $\rho$ with the
 restriction of the Steinberg
Symbol on the $p$-units

$$
E^{(1+i-k)} \times E^{(1- i)} \rightarrow
\mathbb{F}_p(\omega^{2-k})
$$
is non-degenerate.

 We have the following theorem.

\bth Let $p>3$ be a prime and $k$ an even integer, $2\leq k<2p,$
for which hypothesis $H_k$ holds. Then there exists a non-zero
parabolic cohomology class $\psi\in H^1_{par}(SL_2(\Z),
Symm^{k-2}(\F_p^2))$ such that
\begin{enumerate}
\item $\psi|T_q = (1+q^{k-1})\psi$ for $q = 2,3$. \item For
$i=3,5,\ldots,k-3$, we have $L(\psi,i) = \rho(\xi_i)$;
\end{enumerate}
\end{theorem}

The above results were motivated by the joint work of W.\ McCallum
and R.\ Sharifi \cite{MS1},\cite{MS2} and by a well-publicized
conjecture of Sharifi (see \cite{SharifiSlides}.   Indeed,
McCallum and Sharifi \cite{MS2} predict that, assuming Vandiver's
Conjecture, Hypothesis $H_k$ holds whenever $p$ is irregular with
$p|{\frac{B_k}{k}}$.   In this case, Sharifi's conjecture
\cite{SharifiSlides} predicts the truth of Theorem 1.2 with the
Hecke condition (1) strengthened to include all primes $q\ne p$
(not just $q=2,3$). In future work of the author with G.\ Stevens,
we will generalize Theorem 1.2 by constructing a universal
Eisenstein cohomology class $\Psi$ that specializes to the class
$\psi$ of Theorem 1.2. The general statement (1) about Hecke
eigenvalues will follow from the properties of $\Psi$.   Details
will appear later.

 Now, for a set
$S$ of positive primes $q\ne p$, let
$$
H_{k,eis,S}^+ \subseteq H^1_{par}(SL_2(\Z),Symm^{k-2}(\F_p))
$$
denote the subspace of all vectors that are fixed by complex
conjugation and on which the operators $T_q$ for $q\in S$ act with
eigenvalue $1+q^{k-1}$. It follows from Theorem 1.2 that
$H^+_{k,eis,\{2,3\}}$ is positive dimensional whenever hypothesis
$H_k$ holds. The following theorem is a consequence of Theorem
1.2.

\begin{theorem}
Let $(p,k)$ be an irregular pair such that hypothesis $H_k$ holds
and assume $H^+_{k,eis,\{2,3\}}$ is one-dimensional. Then
Sharifi's conjecture is true for the pair $(p,k)$.
\end{theorem}
Finally, we remark that Romyar Sharifi has recently and
independently also proved Theorem 1.1.

I would like to mention that the present work would not have been
possible without the insightful advice and wealth of ideas from my
thesis advisor, Glenn Stevens. I am also extremely grateful for
his immense help with the expository part of the paper and Lemma
7.5. I would also like to thank Romyar Sharifi for his generosity
in sharing copies of transparencies from talks he has given on his
conjectures and also for a number of helpful conversations.

\section{Modular symbols}

Let $\Delta=Div(\mathbb{P}^1(\Q))$ and $\Delta_0\subseteq\Delta$
be the subgroup of divisors of degree $0$.  The group $GL_2(\Q)$
acts by fractional linear transformations on $\Delta$ and
$\Delta_0$. Let
$\Sigma_0(p^n)=\{\left(\begin{smallmatrix}a\;b\\c\;d\end{smallmatrix}\right)\in
M_2^+(\Z)|(a,p)=1,p^n|c\}$ and
$\Sigma_1(p^n)=\{\left(\begin{smallmatrix}a\;b\\c\;d\end{smallmatrix}\right)\in
M_2^+(\Z)|a\equiv 1(mod\;p^n),p^n|c\}$.

 For any right $\Sigma_0(p^n)$-module $M$ we define
a right action of $\Sigma_0(p^n)$ on $\Hom_{\Z}(\Delta_0, M)$ by
$$
\bigl(\phi|\sigma\bigr)(D) = \phi(\sigma D)|\sigma
$$
for all $\sigma \in \Sigma_0(p^n),\ D \in\Delta_0$. The group of
$M$-valued modular symbols over $\Gamma_0(p^n)$ is defined to be
the group
$$
\symbn(M) := \Hom_\Z(\Delta_0,M)^{\Gamma_0(p^n)}.
$$

For each positive integer $m$ we define the Hecke operator
$$
\begin{array}{rl}
T_m : \symbn(M) &\longrightarrow \symbn(M)\\
\phi &\longmapsto \phi|T_m := \sum_i \phi|\delta_i
\end{array}
$$
where the sum is over a complete set of representatives
$\{\delta_i\}_i$ for the left  $\Gamman$-cosets in the double
coset $\Gamman \left(\begin{smallmatrix} 1&0\\ 0&m
\end{smallmatrix}\right)
\Gamman $.

For a cusp $r\in \bP^1(\Q)$, let $\Gamma_r$ be the stabilizer
subgroup in $\Gamman$ of $r$. Then for any  $m\in M^{\Gamma_r}$ we
define $ \wphi_{r,m}: \bP^1(\Q)  \rightarrow M $ by defining
$$
\wphi_{r,m}(s)=\begin{cases} m|\gamma^{-1} &\hbox{\rm if $s=\gamma
r,$ with $\gamma\in\Gamman$} \\
0 & \hbox{\rm otherwise.}
\end{cases}
$$
We extend $\wphi_{r,m}$ by linearity to an additive function
$$
\wphi_{r,m} : \Delta \longrightarrow M
$$
and note that the restriction of $\wphi_{r,m}$ to $\Delta_0$ is an
$M$-valued modular symbol over $\Gamman$ which we denote
$\phi_{r,m}$.

\begin{defn}\label{defnBound}
A modular symbol of the form $\gamma_{r,m}$ will be called a
boundary symbol supported on the $r$-cusps, i.e.\  the cusps that
are $\Gamman$-equivalent to $r$. We define the group of $M$-valued
boundary symbols over $\Gamman$ to be the subgroup
$$
\boundn(M)\subseteq \symbn(M)
$$
generated by the set of all $\phi_{r,m}$, where $r,m$ run over all
pairs with $r\in\bP^1(\Q)$ and $m\in M^{\Gamma_r}$.
\end{defn}

There is also a canonical map $\symbn(M) \longrightarrow
H^1(\Gamman,M)$.   If $\phi\in \symbn(M)$ and $r\in \bP^1(\Q)$,
then the map $\Gamman \longrightarrow M$ defined by $ \gamma
\longmapsto \phi\bigl( (\gamma r) - (r) \bigr) $ is a $1$-cocycle,
whose cohomology class is independent of the choice of $r$. We let
$\pi_\phi$ be that cohomology class.   From the definition it is
clear that for any $r\in\bP^1(\Q)$, the restriction of $\pi_\phi$
to $\Gamma_r$ is trivial. Thus, we have
$$
\pi_\phi \in H^1(\Gamman,M).
$$

We have the following theorem of Ash and Stevens \cite{AS}.

\bth If multiplication by $6$ is invertible on $M$ then there is a
canonical isomorphism $ H^1_c(\Gamma_0(p^n),M) \cong \symbn(M) $.
Moreover, there is a canonical commutative diagram
$$
\begin{array}{ccccccccc}
0&\longrightarrow&\boundn(M) &\longrightarrow &\symbn(M)&\buildrel\pi_{\phi}\over\longrightarrow &H^1_{par}(\Gamman,M)&\longrightarrow&0\\
&&\downarrow &&\downarrow && \Vert &&\\
0&\longrightarrow&H^1_{\partial}(\Gamman,M) &\longrightarrow &
H_c^1(\Gamman,M) &\longrightarrow
&H^1_{par}(\Gamman,M)&\longrightarrow&0
\end{array}
$$
in which the rows are exact, the vertical arrows are isomorphisms,
and all maps commute with the natural action of the Hecke
operators $T_m$ ($m\in \N$).  Here, $H_\partial^1$ is the
``boundary cohomology", which is defined by the exactness of the
second row. \end{theorem}

\section{Manin Symbols}

In the special case where the subgroup $\Gamma_1(p^n)\subseteq
\Gamman$ acts trivially on $M$ we can give a simple description of
$\symbn(M)$ in terms of ``Manin symbols".   We recall that
description in this section.

Let $A$ be a commutative ring.   A group homomorphism
$$
\chi : \Gamman \longrightarrow A^\times
$$
will be called a nebentype character if $\chi$ is trivial on
$\Gamma_1(p^n)$.  Any nebentype character $\chi$ extends uniquely
to a multiplicative map
$$
\chi : \Sigma_0(p^n) \longrightarrow A^\times
$$
that is trivial on $\Sigma_1(p^n)$.

\begin{defn}\label{defnNeben}
Let $A$ be a ring and $M$ be an $A$-module endowed with a right
action of $\Sigman$. We say that $\Sigman$ acts via the nebentype
character $\chi$ if for all $\gamma \in \Sigman$ and all $m\in M$
we have $ m|\gamma = \chi(\gamma)\cdot m. $
\end{defn}

The group $SL_2(\Z)$ acts by right matrix multiplication on the
additive group of row vectors $\left(\Z/p^n\Z\right)^2$ .   The
orbit of $(0,1)$ is the set
$$
X_n := \Xn =  \left\{(x,y)\in (\Z/p^n\Z)^2\,\biggm|\,
(x,y,p)=1\right\}.
$$
The stabilizer of $(0,1)$ is the subgroup $\Gamma_1(p^n)$, which
is a normal subgroup of $\Gamma_0(p^n)$.   Thus $\Gamma_0(p^n)$
also acts on $X_n$ on the left.   In fact, if $\gamma =
\left(\begin{smallmatrix} a&b&\\c&d
\end{smallmatrix}\right)
$ then the left action of $\gamma$ on $X_n$ is given by scalar
multiplication by $d$: $ \gamma \bx \mapsto d\cdot \bx. $ For the
rest of this section $M$ will be an $A$-module on which $\Sigman$
acts via the nebentype character $\chi : \Sigman\longrightarrow
A^\times$.

\begin{defn}\label{defMsymb}
A function
$$
e : X_n \lra M
$$
is called an $M$-valued Manin symbol over $\Gamma_0(p^n)$ if $e$
satisfies the following ``Manin relations" for all $\bx = (x,y)\in
X_n$ and $\lambda \in \left(\Z/p^n\Z\right)^\times$:
\begin{itemize}
\item[(1)] $e(\lambda \bx) = \chi(\lambda)\cdot e(\bx)$;
\item[(2)] $e(x,y) + e(y,-x) = 0$; and \item[(3)] $e(x,y) +
e(y,-x-y) + e(-x-y,x) = 0$.
\end{itemize}
We denote by
$$
\Maninn(M).
$$
the group of all $M$-valued Manin symbols over $\Gamman$.
\end{defn}

\bre \label{even} Two applications of the second Manin condition
show that if $e$ is a Manin symbol, then $e(\bx) = e(-\bx)$ for
every $\bx\in X_n$.   So if $e\ne 0$ then the first condition
implies  $\chi$ must be {\it even}, i.e.\ $\chi(-1) = 1$. \ere

Now fix a section $X_n \longrightarrow SL_2(\Z)$, $\bx \mapsto
\gamma_\bx$, so that
$$
(0,1)\gamma_\bx = \bx
$$
for every $\bx\in X_n$.    Also let $D_\bx \in \Delta_0$ be given
by
$$
D_\bx := \gamma_\bx\cdot \bigl((\infty) - (0)\bigr).
$$
If $M$ is a right $\Gamma_0(p^n)$-module and $\phi \in \symbn(M)$,
then we define $\bm_\phi : X_n \longrightarrow M$ by
$$
\bm_\phi(\bx) := \phi(D_\bx)
$$
and note that this is well-defined independent of our choices of
the $\gamma_\bx$. We have the following reformulation of a theorem
of Manin \cite{Manin1}.

\begin{theorem}
The map $\phi \longmapsto \bm_\phi$ induces an isomorphism
$$
\bm : \symbn(M) \lra \Maninn(M)
$$
for every right $\Gamma_0(p^n)$-module $M$.
\end{theorem}

We use this isomorphism to transfer the action of the Hecke
operators $T_m$ on $\symbn(M)$ to an action on $\Maninn(M)$.   The
following theorem of Merel \cite{Merel} gives an ``explicit"
description of this action.

\begin{theorem}\label{merel}
Let $m$ be a positive integer and let
$$
H_m : = \left\{\, \left(\begin{matrix} a&b\\ c&d
\end{matrix}
\right) \bigg |
\begin{array}{l}
a,b,c,d\in\Z\\
a>b\ge 0,\ d>c\ge 0\\
ad-bc = m
\end{array}
\, \right\}.
$$
Then for every $e\in \Maninn(M)$ we have
$$
(e|T_m)(\bx) = \sum_{\delta\in H_m} e(\bx\delta).
$$
\end{theorem}

For future reference we record the following corollary.

\begin{corollary}\label{merelcor}
For arbitrary $e\in\Maninn(M)$ we have
$$\begin{array}{rcl}
(e|T_2)(\bx) &=& e(x,2y) + e(2x,y) + e(x+y,2y) + e(2x,x+y)\\
(e|T_3)(\bx) &=& e(x,3y) + e(3x,y) + e(x+y,3y) + e(3x,x+y)\\
&& + e(x-y,3y)  + e(3x,x-y).
\end{array}
$$
\end{corollary}

\begin{proof}  We easily verify that
$$
H_2 = \left\{ \psmat 1&0\\ 0 & 2 \epsmat, \psmat 2&0\\ 0 & 1
\epsmat, \psmat 1&0\\ 1 & 2 \epsmat, \psmat 2&1\\ 0 & 1 \epsmat
\right\} \quad \hbox{\rm and} \quad H_3 = \left\{ \psmat 1&0\\ 0 &
3 \epsmat, \psmat 3&0\\ 0 & 1 \epsmat, \psmat 1&0\\ 1 & 3 \epsmat,
\psmat 1&0\\ 2 & 3 \epsmat, \psmat 3&1\\ 0 & 1 \epsmat, \psmat
3&2\\ 0 & 1 \epsmat, \psmat 2&1\\ 1 & 2 \epsmat \right\}.
$$
The description of $T_2$ is then an immediate consequence of
Merel's theorem. On the other hand, Merel's theorem implies
$$\begin{array}{rcl}
(e|T_3)(\bx) &=& e(x,3y) + e(3x,y) + e(x+y,3y) + e(3x,x+y) \\
&& + e(x+2y,3y) + e(3x,2x+y) + e(2x+y,2y+x).
\end{array}
$$
Now consider the following matrix:
$$\pmat
e(x+2y,3y) & e(x-y,x+2y) & 0  \\
e(3x,2x+y) & 0  & e(2x+y,x-y)  \\
e(2x+y,x+2y) & e(x+2y,x-y) & e(x-y,2x+y) \epmat
$$
From the Manin relations, we see that the sum of the three rows
are $e(x-y,3y)$, $e(3x,x-y)$, and $0$, respectively.   Thus the
sum of all  the entries of the matrix is $e(x-y,3y)+e(3x,x-y)$. On
the other hand, the second and third columns sum to $0$. Hence the
sum of all the entries of the matrix is the sum of the first
column.  We therefore have
$$
e(x+2y,3y) + e(3x,2x+y) + e(2x+y,x+2y) = e(x-y,3y) + e(3x,x-y).
$$
Substituting this into the above expression for $e|T_3$ gives us
$$\begin{array}{rcl}
(e|T_3)(\bx) &=& e(x,3y) + e(3x,y) + e(x+y,3y) + e(3x,x+y) \\
&& + e(x-y,3y) + e(3x,x-y) ,
\end{array}
$$
proving our claim for $e|T_3$.
\end{proof}

We conclude this section with a discussion of  boundary symbols
supported on the $\infty$-cusps (see Definition \ref{defnBound}).

\begin{defn}
We say that a Manin symbol $e\in \Maninn(M)$ is supported at
$\infty$ if $e$ satisfies the following condition:
$$
e(\bx) = 0 \quad \hbox{\rm for all $\bx = (x,y)\in X_n$ with
$xy\ne 0$}.
$$
We denote by
$$
\Manininf(M) \subseteq \Maninn(M)
$$
the subgroup of all Manin symbols supported at $\infty$.
\end{defn}

It is easy to describe the action of the Hecke operators on
$\Manininf(M)$. The result is the following.

\begin{proposition}\label{heckeBound}
The subgroup $\Manininf(M)$ is an eigen-submodule for the action
of the Hecke operators $T_m$  on $\Maninn(M)$.   Moreover, we
have:
\begin{itemize}
\item if $\ell\ne p$ is prime, then the eigenvalue of $T_\ell$ is
$\ell+\chi(\ell)$; \item the eigenvalue of $T_p$ is $p$.
\end{itemize}
\end{proposition}

The proof is an easy computation from the definitions.   We note
that the formal Dirichlet series of Hecke acting on $\Manininf(M)$
is given by
$$
\sum_{m=1}^\infty T_m m^{-s} =
\left(1-p^{1-s}\right)^{-1}\prod_{\ell\ne
p}\left(1-\ell^{1-s}\right)^{-1}\left(1-\chi(\ell)\ell^{-s}\right)^{-1},
$$
which we may write suggestively in the form
$$
\sum_{m=1}^\infty T_m m^{-s} = \zeta(s-1)\cdot L(s,\chi).
$$
The right hand side is the $L$-function of an Eisenstein series
$E_\chi$ defined over the ring $A$.

\section{$K_2$ of cyclotomic integer rings}\label{kthry}

First, we recall the definition of Milnor's $K_2$-group of a
commutative ring $R$.

\begin{defn}   The second Milnor $K$-group of $R$
is defined as
$$
K_2^M(R) := (R^\times\otimes_{\Z}R^\times)/I_2,
$$
where $I_2$ is the subgroup of $R^\times \otimes_\Z R^\times$
generated by the set
$$
\left\{a_1\otimes a_2\in
R^\times\otimes_{\Z}R^\times\,\biggm|\,a_1+a_2 \in \{0,1\}
\right\}.
$$
The Steinberg Symbol is defined as the canonical map
$$
\{\ ,\ \} : R^\times \times R^\times \longrightarrow K_2^M(R).
$$
\end{defn}

We write the multiplication in $K_2^M(R)$ additively and note that
for $a,b\in R^\times$ we have
$$
\{a,b\} + \{b,a\} = 0
$$
as an easy consequences of the relation $\{ab,-ab\} = 0$.   This
skew symmetry of the Steinberg symbol will be used throughout the
paper.

Denote by $\mu_n$ the group of $p^n$-th roots of unity in
$\C^\times$ and fix $\zeta_n\in \mu_n$ a primitive $p^n$-th root
of unity.  Let $K_n :=\Q(\mu_{n})$, $R_n=\Z[\mu_{n},\frac{1}{p}]$,
and $G_n := \hbox{\rm Gal}(K_n/\Q)$.   We associate to each $a\in
(\Z/p^n\Z)^\times$ the element $\sigma_a \in G_n$ for which
$\zeta_{n}^{\sigma_a}=\zeta^a$. As in the introduction, we let
$\wK_2(R_n) = K_2^M(R_n)/(2\hbox{\rm -torsion})$. Note that
$\wK_2(R_n)$ has a natural action of the Galois group $G_n$.

In the next section, we will need the following lemma.

\begin{lemma}\label{trivunits}
For any $x,y\in \Z/p^n\Z$ with $x\ne 0$ we have
$$
\{1-\zeta_n^x,\zeta_n^y\} = 0.
$$
\end{lemma}

\begin{proof}
If $x\in(\Z/p^n\Z)^\times,$ then we may choose $a\in \Z$ such that
$ax = 1$ in $(\Z/p^n\Z)$. But  then $\{1-\zeta_{n}^x,\zeta_{n}^y\}
= \{1-\zeta_{n}^x,\zeta_{n}^{axy}\} = ay\cdot \{1-\zeta_{n}^x,
\zeta_{n}^x\}=0$, proving the lemma in this special case. In the
general case we may write $x \equiv p^{k}u$ (mod $p^n$) for
$u,k\in \Z$ with $p\hbox{$\not|$}u$ and $0 \le k < n$. Then,
$$1-\zeta_{n}^x = \prod_{\alpha\in(\Z/p^n\Z)^\times\atop \alpha\equiv
1\ (p^{n-k})}1-\zeta_{n}^{u\alpha},$$ and the relation will now
follow from the special case.
\end{proof}

\section{A $\wK_2(R_n)$-valued Manin symbol}\label{K-Manin}

Define the {\emph Artin nebentype character} to be the character
$$
\begin{array}{rl}
\chi : \Sigman & \longrightarrow \Z[G_n]^\times\\
\gamma &\longmapsto \sigma_a
\end{array}
$$
where, as always, $a$ is the upper left corner of $\gamma$. We let
$\Sigman$ act on $\wK_2(R_n)$ via the Artin nebentype $\chi$.

\begin{theorem} \label{thmKManin}
The function $e_n : X_n \longrightarrow \wK_2(R_n)$ defined by
$$
e_n(x,y) =
\begin{cases}
\{1-\zeta_n^x,1-\zeta_n^y\} &\hbox{\rm if $x,y \neq 0,$}\\
 0 & \hbox{\rm otherwise,}
\end{cases}
$$
is a Manin symbol: $e_n\in \Maninn(\wK_2(R_n))$. Moreover, for
$q=2,3$
$$
e_n\biggm|\biggl(T_q - \bigl(q+\chi(q)\bigr)\biggr) \in
\Manininf(\wK_2(R_n)) .
$$
\end{theorem}

\begin{proof}   We begin by verifying that $e_n$ satisfies the three Manin conditions (see Definition \ref{defMsymb}).   Let $\bx = (x,y) \in X_n$ be fixed and note that the first two Manin conditions are trivially satisfied if $xy=0$.

For $\lambda\in \left(\Z/p^n\Z\right)^\times$ we have
$$\begin{array}{rl}
e_n(\lambda\bx) &=  \{1-\zeta_n^{\lambda x},1-\zeta_n^{\lambda
y}\} =
\{1-\zeta_n^x,1-\zeta_n^y\}^{\sigma_\lambda}\\
&=\chi(\lambda)\cdot e_n(\bx)
\end{array}
$$
So the first Manin condition is satisfied.

From the properties of the Steinberg symbol we have, for $xy\ne
0$,
$$
\begin{array}{rcl}
e_n(x,y) & \sim & \{\zeta_n^x-1,1-\zeta_n^y\} = -\{1-\zeta_n^y,\zeta_n^x - 1 \} \\
&=& -\{1- \zeta_n^y, \zeta_n^x\}\cdot \{1-\zeta_n^y,\ 1-\zeta_n^{-x}\}\\
&=& -\{1-\zeta_n^y,1-\zeta_n^{-x}\}\\
&=& -e_n(y,-x).
\end{array}
$$
Here $\sim$ denotes congruence modulo $2$-torsion in $K_2^M(R_n)$
(recall that $\wK_2(R_n) = K_2^M(R_n)/\hbox{\rm ($2$-torsion)}$).
We used Lemma \ref{trivunits} to derive the second to last
equality. This proves the second Manin condition.

To verify the third Manin condition,
$$
e_n(x,y) + e_n(y,-x-y) + e_n(-x-y,x) = 0
$$
we consider cases. If $x= 0$ then $y\ne 0$, and in that case
$e_n(y,-y) = -e_n(y,y)$ by the second Manin condition. But by the
skew symmetry of the Steinberg symbol we have $e_n(y,y)\sim 0$, so
$e_n(y,-y) \sim 0$ and the third Manin condition is satisfied in
this case.  Similarly, it is satisfied if either $y=0$ or $x+y=0$.
So we may assume $x,y,$ and $x+y$ are all non-zero..   In this
case, we have the identity
$$
\frac{\zeta_{n}^y(1-\zeta_{n}^x)}{1-\zeta_{n}^{x+y}}+\frac{1-\zeta_{n}^y}{1-\zeta_{n}^{x+y}}=1.
$$
From the Steinberg relations we then have
$$
\left\{\frac{\zeta_{n}^y(1-\zeta_{n}^x)}{1-\zeta_{n}^{x+y}},
\frac{1-\zeta_{n}^y}{1-\zeta_{n}^{x+y}}\right\} =0.
$$
Bimultiplicativity of the Steinberg symbol, Lemma \ref{trivunits},
and the skew symmetry of $e_n$ imply
$$
e_{n}(x,y)-e_{n}(x+y,y) -e_{n}(x,x+y) = 0.
$$
Now apply the second Manin condition to the last two terms to
obtain
$$
e_{n}(x,y)+e_{n}(y,-x-y)+e_{n}(-x-y,x)= 0,
$$
and the third Manin condition is proved.   This proves $e_n$ is a
Manin symbol.

To compute the Hecke operators, we use the corollary
\ref{merelcor} to Merel's theorem \ref{merel}.

\ble\label{primHecke} Let $q=2$ or $q=3$.   Then for all $\bx =
(x,y)$ with $xy\ne 0$ we have
$$
(e_n|T_q)(\bx) = (q + \chi(q))\cdot e_n(\bx).
$$
\ele

\begin{proof}

For $q=2$ we have, from Corollary \ref{merelcor}
$$
(e_n|T_2)(x,y) = e_n(x,2y) + e_n(2x,y) + e_n(x+y,2y) +
e_n(2x,x+y).
$$

If $x+y=0$ this says $ (e_n|T_2)(x,-x) = e_n(x,-2x) + e_n(2x,-x) $
which vanishes by skew symmetry and the second Manin condition.
But also $e_n(x,-x) = 0$, so we have
$$
(e_n|T_2)(x,-x) = (2+\chi(2))\cdot e_n(x,-x)
$$
since both sides of this equation vanish.

If $x+y\ne 0$ then we use the identity
$$
\frac{(1-\zeta_{n}^{x+y})(1-\zeta_{n}^x)}{1-\zeta_{n}^{2x}}+\zeta_{n}^x\frac{(1-\zeta_{n}^{2y})(1-\zeta_{n}^x)}{(1-\zeta_{n}^{2x})
(1-\zeta_{n}^y)}=1.
$$
This implies $e_{n}(x,2y)+e_{n}(2x,y)+e_{n}(x+y, 2y) -
e_{n}(x+y,2x) = e_{n}(x,y) +  e_{n}(2x,2y)  + e_{n}(x+y,y) -
e_{n}(x+y,x)  + e_{n}(x,2x)  + e_{n}(2x,x)$.  According to the
above, the left hand side of this equality is $(e_n|T_2)(x,y)$. So
we have
$$
\begin{array}{rcl}
(e_{n}|T_2)(x,y) &=& 2e_{n}(x,y) + e_{n}(2x,2y) \\
&& -\bigl(e_n(x,y) + e_{n}(y,x+y) + e_{n}(x+y,x)\bigr)\\
&& + \bigl(e_{n}(x,2x)  + e_{n}(2x,x)\bigr)
\end{array}
$$
But the last two lines of the right side of this equation vanish
by the third and second Manin conditions, so we have
$$
(e_n|T_2)(\bx) = (2+\chi(2))\cdot e_n(\bx)
$$
and the assertion for $T_2$ is proved.

For $q=3$ we again use Corollary \ref{merelcor} to obtain
$$
\begin{array}{rcl}
(e_n|T_3)(x,y) &=& e_n(x,3y) + e_n(3x,y) + e_n(x+y,3y)\\
&& + e_n(x-y,3y) + e_n(3x,x+y)  + e_n(3x,x-y).
\end{array}
$$

If either $x+y =0$ then the right hand side simplifies to
$e_n(x,-3x) + e_n(3x,-x) + e_n(-x,-3x) + e_n(3x,x)$ which vanishes
by the skew symmetry of $e_n$. A similar calculation shows that
the right hand side vanishes when $x-y = 0$. Thus in either case,
we have
$$
(e_n|T_3)(x,\pm x) = (3+\chi(3))\cdot e_n(x,\pm x)
$$
since both sides vanish.

So we may assume $x+y, x-y \neq 0$.  In that case, we have the
identity
$$
\zeta_{n}^{y-x}(\zeta_{n}^{2x}+\zeta_{n}^x+1)+(1-\zeta_{n}^{y-x})(1-\zeta_{n}^{x+y})=\frac{1-\zeta_{n}^{3y}}{1-\zeta_{n}^y},
$$
which may be rewritten as
$$
\frac{\zeta_{n}^{y-x}(1-\zeta_{n}^{3x})(1-\zeta_{n}^{y})}{(1-\zeta_{n}^x)(1-\zeta_{n}^{3y})}+
\frac{(1-\zeta_{n}^{y-x})(1-\zeta_{n}^{y})(1-\zeta_{n}^{x+y})}{1-\zeta_{n}^{3y}}=1.$$

This implies
$e_{n}(x,3y)+e_{n}(3x,y)-e_{n}(3y,x+y)-e_{n}(3y,y-x)+e_{n}(3x,x+y)
+e_{n}(3x,y-x) = e_{n}(3x,3y)-e_{n}(y,y-x)-e_{n}(y,x+y) +
e_{n}(x,y-x) + e_{n}(x,y) + e_{n}(x,x+y)$.

The left hand side of this equality is $(e_n|T_3)(x,y)$.   So we
have
$$
\begin{array}{rcl}
(e_{n}|T_3)(x,y) &=&3e_{n}(x,y) + e_{n}(3x,3y) \\
&& + e_n(y,x) + e_{n}(x,y-x) + e_{n}(y-x,y)\\
&& + e_n(y,x) + e_{n}(x,x+y) + e_{n}(x+y,y).
\end{array}
$$
Using the third Manin relation, we see that the bottom two rows of
the right hand side vanish.   Hence
$$
(e_n|T_3)(\bx) = (3 + \chi(3))\cdot e_n(\bx)
$$
and the lemma is proved.
\end{proof}

We mention that the relations on $p$-units we used for the
computations of the Hecke operators $T_2$ and $T_3$ were remarked
by W.McCallum and R. Sharifi \cite{MS1} \cite{SharifiComps}.

We now return to the proof of the theorem.   Let $q=2$ or $q=3$.
It follows from the lemma that the Manin symbol $e :=
e_n|(T_q-(q+\chi(q))$ vanishes on all $\bx = (x,y)\in X_n$ with
$xy\ne 0$. Thus from Proposition \ref{heckeBound} we see that
$e\in \Manininf(\wK_2(R_n))$ and theorem is proved.
\end{proof}

\section{The parabolic cohomology class $\varphi_n$}

For each positive integer $n$ we let
$$
B_n := \boundn(\wKn),\ \ {\rm and}\ \ S_n := \symbn(\wKn)
$$
and consider the exact sequence
$$
0 \longrightarrow B_n \longrightarrow S_n
\buildrel\pi\over\longrightarrow H^1_{par}(\Gamman,\wKn)
\longrightarrow 0.
$$
Let $\phi_n \in S_n$ be the modular symbol associated to the Manin
symbol $e_n$ defined in the previous section. We define
$$
\varphi_n := \pi_{\phi_n} \in H^1_{par}(\Gamman,\wKn)
$$
to be the image of $\phi_n$.   Theorem 1.1 in the Introduction
follows now immediately from Theorem \ref{thmKManin}.

We raise a number of key questions.

\noindent {\bf Questions:}
\begin{itemize}
\item[(1)]  Is $\varphi_n$ an eigenclass for all the Hecke
operators? \item[(2)]  Under what conditions can we say
$\varphi_n\ne 0$? \item[(3)]  What relations exist among the
$\varphi_n$ as $n$ varies? \item[(4)]  Assuming $\varphi_n\ne 0$,
is there a Hecke eigensymbol $\psi_n\in S_n$ lifting $\varphi_n$?
\end{itemize}

In fact, in future work, we will draw a connection between the
modular symbols $\phi_n$ and the Eisenstein distribution \cite{St}
and use this connection to show that $\varphi_n$ is an eigenclass
satisfying
$$
\varphi_n|T_q = (q + \chi(q))\cdot \varphi_n
$$
for all primes $q$, where we understand that $\chi(q) = 0$ when
$q=p$. Thus (1) has an affirmative answer.

The answers to questions (2) to (4) are closely tied to some very
beautiful recent conjectures of Romyar Sharifi
\cite{SharifiEisIdeal}, which in turn are motivated by work of
Ohta (see \cite{Ohta}). Sharifi's ideas suggest that (2) is
closely connected to the structure of the class group of $K_n :=
\Q(\zeta_n)$.   The answer to (3) should be given (for $m\le n$)
in terms of the transfer map $\wKn \longrightarrow \wK_2(R_m)$ on
$K$-theory composed with corestriction of the cohomology of
$\Gamma_0(p^n)$ to $\Gamma_0(p^m)$.   Finally, we expect the
answer to (4) to be negative, which corresponds to an expectation
that there should be lots of fusion between the boundary
cohomology and the parabolic cohomology.

\section{Special values of $L$-functions}

In this section, following \cite{AS}, we define the universal
$L$-value of a modular symbol and describe a few of its
properties.   In particular we will see that the Manin symbols are
universal $L$-values.

Let $R$ be a ring and $M$ be an $R$-module endowed with a right
action of $SL_2(\Z)$.

\begin{defn} (Universal $L$-values)   Let $\phi\in H^1_c(SL_2(\Z),M)$ be
a modular symbol.  We define $\Lambda(\phi)\in M$ by:
$$
\Lambda(\phi) := \phi\bigl((\infty) - (0)\bigr)
$$
and call $\Lambda(\phi)$ the universal $L$-value of $\phi$.
\end{defn}

We define $M^\ast$ to be the $R$-dual of $M$:  $M^\ast :=
\Hom_R(M,R)$ with $SL_2(\Z)$ acting on the right as: $
(\lambda|\sigma)(m) = \lambda(m|\sigma^\prime) $ where
$\sigma\mapsto \sigma^\prime$ is the adjugate involution $\psmat
a&b\\c&d\epsmat \mapsto \psmat d&-b\\-c&a\epsmat$. For example,
let
$$
W_r(R) = \left\{\,F\in R[X,Y]\,\biggm|\,\hbox{$F$ is homogeneous
of degree $r$}\right\}
$$
with $M_2^+(\Z)$ acting by the formula $ (F|\sigma)(X,Y) =
F((X,Y)\sigma^\prime) $ and let
$$
V_r(R) := W_r(R)^\ast.
$$

As in [AS86] we make the following definition.

\begin{defn}\label{Lvals}
Let $\phi\in H^1_c(SL_2(\Z),V_r(R))$.   We define the special
$L$-values $L(\phi,i+1)\in R$ for $i=0,1,\ldots,r$ by
$$
L(\phi, i + 1) := \left\langle
\Lambda(\phi),(-1)^iX^{r-i}Y^i\right\rangle
$$
where $\langle\ ,\ \rangle : V_r(R) \times W_r(R) \lra R$ is the
canonical pairing.
\end{defn}

If $r!$ is invertible in $R$ then there is a unique
$M_2^+(\Z)$-equivariant perfect pairing
$$
\langle\,,\,\rangle : W_r(R) \times W_r(R) \lra R
$$
with respect to which
$$
\left \langle \psmat r \\ i \epsmat
X^iY^{r-i},(-1)^jX^{r-j}Y^j\right \rangle =
\begin{cases}
1 & \hbox{\rm if $i=j$},\\
0 & \hbox{\rm otherwise}.
\end{cases}
$$
Thus, when $r!$ is invertible in $R$ we have $V_r(R) \cong W_r(R)$
and we may regard $\Lambda(\phi)$ as an element of $W_r(R)$.  With
the above identifications we then have
$$
\Lambda(\phi) := \sum_{i=0}^r \psmat r \\ i \epsmat
L(\phi,i+1)X^iY^{r-i},
$$
consistent with the conventions of [AS86].  Note however that
definition \ref{Lvals} is meaningful for any commutative ring $R$
-- we do not need to assume $r!$ is invertible in $R$.

Now consider the general case.   Let $\Gamma$ be any congruence
subgroup of $SL_2(\Z)$, let $M$ be a $\Gamma$-module,  and let
$\varphi \in \symb_\Gamma(M)$ be an $M$-valued modular symbol over
$\Gamma$. To define the universal $L$-value of $\varphi$ we first
induce to $SL_2(\Z)$ using Shapiro's Lemma and then take the
universal $L$-value of the induced modular symbol.

More precisely, we define the induced module of $M$ to be the
module
$$
I(M) := \left\{\,f : SL_2(\Z)\to M\ \biggm|\,f(\gamma x) =
f(x)|\gamma^{-1},\ \forall \gamma\in \Gamma,\ x\in
SL_2(\Z)\,\right\}
$$
with $SL_2(\Z)$ acting by $(f|g)(x) = f(xg^{-1})$. The Shapiro
isomorphism gives us a canonical isomorphism
$$
H^1_c(\Gamma,M) \buildrel I \over \lra H^1_c(SL_2(\Z),I(M))
$$
which is given explicitly on modular symbols by $ I :\varphi
\longmapsto I(\varphi)$ where $I(\varphi) : \Delta_0\lra I(M)$ is
given by
$$
I(\varphi)(D)(x) = \varphi(xD)
$$
for $D\in \Delta_0$ and $x\in SL_2(\Z)$.

\begin{defn}  Let $\varphi\in H^1_c(\Gamma,M)$.  Then the universal $L$-value
of $\varphi$ is defined to be
$$
\Lambda(\varphi) := I(\varphi)\bigl((\infty) - (0)\bigr).
$$
In other words we set $\Lambda(\varphi) = \Lambda(I(\varphi))$.
\end{defn}

In the special case where $M$ is a $\Gamman$-module on which
$\Gamman$ acts via a nebentype character $\chi$ we may identify
$I(M)$ with the module of functions $f: X_n \lra M$ satisfying
$f(d\bx) = \chi(d)\cdot f(\bx)$.   Thus we have a natural
inclusion
$$
\Maninn(M) \hookrightarrow I(M).
$$
In fact, we have the following simple proposition, whose proof we
leave to the reader.

\begin{proposition}\label{LManin}
Let $M$ be an $R$-module on which $\Sigman$ acts via a nebentype
character $\chi$ and let $\phi\in H^1_c(\Gamman,M)$ be a modular
symbol.   Then with the above identifications, we have
$$
\Lambda(\phi) = e_\phi.
$$
\end{proposition}

Finally, we turn to the problem of defining special $L$-values of
parabolic cohomology classes. For this we need to understand the
module $V_r(R)^{\Gamma_\infty}$.   For the rest of this section
$R$ will be a ring of characteristic $p$ and $r\ge 0$ will an even
integer.   Let $\{\lambda_i\}_{i=0}^r$ in $V_r(R)$ be the dual
basis to $\{(-1)^iX^{r-i}Y^i\}_{i=0}^r$ in $W_r(R)$. We leave the
simple proof of the following lemma to the reader.

\begin{lemma}\label{invlemma}
Let $r$ be an even integer.   If $r < p$ then $\lambda_r$ spans
$V_r(R)^{\Gamma_\infty}$.   If $p \le r < 2p$ then $\lambda_r$ and
$\lambda_{p-1}$ span $V_r(R)^{\Gamma_\infty}$.
\end{lemma}

\begin{proposition}\label{boundvals}
Let $\Gamma = SL_2(\Z)$, let $i,r$ be positive integers with $r$
even and $0 \le  i \le r$.  Let
$\varphi\in\bound_{\Gamma}(V_r(R)))$ be a boundary symbol. Then
\begin{itemize}
\item[(a)] If $r < 2p$ then $ L(\varphi,i+1) = 0\ \ \hbox{\rm for
$i\not\equiv 0, r$ (mod $p-1$).} $ \item[(b)] For arbitrary $r$,
if $\varphi|T_p = \varphi$ then $ L(\varphi,i+1) = 0\ \ \hbox{\rm
for $i\ne  0, r$.} $
\end{itemize}

\end{proposition}

\begin{proof}  Fix $i,r$ as in the statement of the proposition and let
$\varphi\in \bound_{\Gamma}(V_r(R)))$. Since $SL_2(\Z)$ has only
one cusp, we have $\varphi = \varphi_{\infty,\lambda}$ for some
$\lambda\in V_r^{\Gamma_\infty}$. By lemma \ref{invlemma}  there
are constants $a,b\in R$ such that $\lambda = a\lambda_r$, if $r <
p$, and $\lambda = a\lambda_r + b\lambda_{p-1}$ if $p < r < 2p$.
But then
$$
\Lambda(\varphi) = \lambda  - \lambda | \psmat 0 & -1\\ 1 & 0
\epsmat =
\begin{cases}
a (\lambda_r -  \lambda_0) &\hbox{ if $r < p$}\\
a (\lambda_r -  \lambda_0) + b(\lambda_{p-1} -  \lambda_{r-p+1})
&\hbox{ if $p < r < 2p$}.
\end{cases}
$$
For $0 \le  i  \le  r$ we have $L(\varphi,i+1)$ is the coefficient
of $\lambda_i$, and by inspection we have, for $r< 2p$,
$$
L(\varphi, i+1) =  0
$$
for $i\not\equiv 0,r$ (mod $p-1$).   This proves (a).

To prove (b) we write $\lambda = \sum_{i=0}^r a_i\lambda_i$ with
coefficients in $R$ and let $m$ be the smallest index for which
$a_m \ne 0$.  Then, since $\varphi|T_p = \varphi$ we have
$$
\lambda = (\varphi|T_p)(\infty) = \left(\varphi|\psmat p&0\\ 0&1
\epsmat \right)(\infty) + \sum_{k=0}^{p-1}\left(\varphi|\psmat
1&k\\ 0&p \epsmat \right)(\infty) = \lambda \biggm| \left(\psmat
p&0\\ 0&1 \epsmat  + \sum_{k=0}^{p-1} \psmat 1&k\\ 0&p \epsmat
\right).
$$
Comparing the coefficients of $\lambda_m$ on both sides of this
equation we obtain $p^{m+1} + p^{r-m} = 1$ in the ring $R$.   But
since $p=0$ in $R$ this can only happen if $r=m$.
\end{proof}

\begin{defn}
For $\psi\in H_{par}^1(\Gamma,V_r(A))$, with $r < 2p$, we choose
$\wpsi \in \symb_\Gamma(V_r(A))$ to be an arbitrary lift of $\psi$
and define
$$
L(\psi,i+1) = L(\wpsi,i+1) \ \ \hbox{\rm for $0 \le  i \le r$ with
$i \not\equiv 0,r$ (mod $p-1$)}.
$$
We call these the special $L$-values of $\psi$.
\end{defn}

By (a) of the last proposition, these special $L$-values are
well-defined, independent of the choice of $\wpsi$.  In fact, by
(b) of the proposition, if $\wpsi|T_p = \wpsi$ then {\it all} of
the $L$-values $L(\psi,i+1)$ are well-defined in the range $0< i <
r$.

\section{Proof of Theorem 1.2}

In this section we prove Theorem 1.2 of the introduction. We take
$n=1$ and suppress the subscript $n$ from the notation. Thus
$\zeta = \zeta_1$ is a primitive $p$th root of unity, $R = R_1 =
\Z\left[\zeta,{1\over p}\right]$, $G= G_1$ is the galois group of
$\Q(\zeta)/\Q$,  $\phi = \phi_1$, $\varphi=\varphi_1$, $e = e_1$,
and $X = X_1 = (\F_p^2)^\prime$.

We also let $k\ge 2$ be an even  integer and set $g = 2 - k$. The
semigroup $M_2^+(\Z)$ acts on $V_{k-2}$ and $I_{r}$ (for any $r$)
by
$$
\begin{array}{rcl}
(F|\sigma)(X,Y) & := &   F((X,Y)\sigma^\prime)\\
(f|\sigma)(\bx)   & := &   f(\bx\sigma^\prime).
\end{array}
$$
for $F\in V_{k-2}$, $f\in I_r$, and $\sigma \in M_2^+(\Z)$, where
in the latter case we take $f(\bx\sigma^\prime) = 0$ in case
$\bx\sigma^\prime = 0$.

There is a natural map of $M_2^+(\Z)$-modules: $W_{k-2}  \lra
I_{k-2}$ defined by sending a polynomial to the function it
represents.    In \cite{AS} it is shown that this map is injective
if $k-2 < p$ and is surjective otherwise.   By duality we obtain a
natural map
$$
{I_g} \buildrel \beta \over \lra V_{k-2}(g),
$$
which is surjective if $k-2<p$ and injective otherwise. Here we
are using the simple fact that
$$
I_g \cong I_{k-2}^\ast(g).
$$
Indeed, the pairing
$$
\begin{array}{rl}
I_{k-2} \times I_g &\lra \F_p\\
(f_1,f_2) &\longmapsto \displaystyle \sum_{\bx \in X} f_1(\bx)
f_2(\bx)
\end{array}
$$
can be seen to induce the above isomorphism. Here, we are using
the notation $M(g) := M \otimes {\det}^g$ for any
$M_2^+(\Z)$-module $M$.  See Lemma 3.2 of \cite{AS} for more
details.

We now turn to the proof of Theorem 1.2.   So we assume Vandiver's
conjecture for $p$, and suppose $2 \le k\le 2p$ and that
hypothesis $H_k$ from the introduction is satisfied.   Thus we
have a $G$-equivariant map
$$\rho:  K_2(R) \rightarrow
\F_p(\omega^{g}).
$$
and an odd integer $i$ with $1 < i < k-1$  such that $\rho(\xi_i)
\ne 0$.   Here $\xi_i := \{\eta_{k-i},\eta_i\}$ as in the
introduction.

We let $\Gamma := SL_2(\Z)$ and define the map
$$
\lambda :  H^1_c(\Gamma_0,\wK_2(R))\longrightarrow
H^1_{par}(\Gamma,V_{k-2})
$$
to be the composition
$$
\lambda : H^1_c(\Gamma_0,\wK_2(R))\buildrel Sh_\rho\over
\longrightarrow H^1_c(\Gamma, I_{g}) \buildrel \nu \over
\longrightarrow H^1_{c}(\Gamma,V_{k-2})
$$
where $Sh_\rho$ is the map induced by $\rho$ and the Shapiro
isomorphism, and $\nu$ is the composition of $\beta$ and the
``twist map"
$$
H_c^1(\Gamma,V_{k-2}(g))\buildrel \tau \over \lra
H^1_{par}(\Gamma,V_{k-2})
$$
which is defined as the identity map on the underlying cohomology
groups. Note, however, that $\tau$ does not commute with the
action of the Hecke operators $T_m$. Indeed, we have $
\tau(\varphi|T_m) = m^g\cdot \tau(\varphi)|T_m $ for every
$\varphi\in H_c^1(\Gamma,V_{k-2}(g))$.

Now let $\phi$ be the modular symbol defined in section
\ref{K-Manin} and set
$$
\wpsi := \lambda(\phi) \in H^1_{c}(\Gamma,V_{k-2}) \quad{\rm and}
\quad  \psi = \hbox{\rm the image of $\wpsi$ in
$H^1_{par}(\Gamma,V_{k-2})$}.
$$
It is proved in \cite{AS}, that $Sh_\rho$ is Hecke equivariant, so
from the last paragraph we have
$$
m^g\cdot \wpsi|T_m = \lambda(\phi|T_m)
$$
for every $m\in \N$.  In particular, for $q=2,3$ we have
$$\begin{array}{rl}
\psi|T_q &= q^{k-2}\lambda\left((q+\chi(q))\cdot \phi\right) = q^{k-2}(q+ q^{g})\cdot \psi\\
&= (1+q^{k-1})\psi
\end{array}
$$
proving (1) of Theorem 1.2.

To compute the special $L$-values $L(\psi,i+1)$ for $2\le i < k-2$
we let $\varphi := Sh_\rho(\phi)$ and use Proposition \ref{LManin}
to conclude that $ \Lambda(\varphi) = \rho\circ e. $ Thus
$\Lambda(\varphi)$ is the function
$$
\begin{array}{rcccl}
\Lambda(\varphi) &:& X &\lra& \F_p\\
&&(x,y) & \longmapsto & \rho\bigl(\{1-\zeta^x,1-\zeta^y\}\bigr) \
\hbox{if $xy\ne0$}.
\end{array}
$$
Now let $\wpsi : = \beta(\varphi)$.    Then $\Lambda(\wpsi) =
\beta\left(\Lambda(\varphi))\right)$.   Thus, we have
$$
\Lambda(\wpsi)=  \sum_{i=0}^{k-2}\lambda_i\sum_{(x,y)\in
X}(-1)^ix^{k-2-i}y^{i}\rho(e(x,y))
$$
which implies
$$
\begin{array}{rl}
(-1)^iL(\wpsi,i+1) & := \displaystyle \sum_{(x,y) \in (\F_p^\times)^2} y^ix^{k-2-i}\rho\bigl(\{1-\zeta^x,1-\zeta^y\}\bigr)\\
&= \displaystyle
\rho\left(\left\{\prod_{x\in\F_p^\times}(1-\zeta^x)^{x^{k-2-i}},
\prod_{y\in\F_p^\times}(1-\zeta^y)^{y^{i}}\right\}\right)\\
&= \displaystyle \rho\left(\left\{\prod_{\sigma\in
G}(1-\zeta)^{\omega^{k-2-i}(\sigma)\sigma}, \prod_{\sigma\in
G}(1-\zeta)^{\omega^{i}(\sigma)\sigma}\right\}\right)
\end{array}
$$
Finally, we recall that $\xi_j \in K_2(R)$ was defined as $\xi_j
:= \{\eta_{k-j},\eta_{j}\}$ where $\eta_j\in E_j$ is the
projection of $(1-\zeta)$ to $E_j := E^{(1-j)}$.    Recalling that
the idempotent projecting to $E^{(1-j)}$ is $ {1\over
p-1}\sum_{\sigma\in G} \omega^{j-1}(\sigma)\sigma, $ we conclude
that
$$
\eta_j^{-1} = \prod_{\sigma\in
G}(1-\zeta)^{\omega^{j-1}(\sigma)\sigma}.
$$
Thus
$$
(-1)^iL(\wpsi,i+1) =
\rho\left(\left\{\eta_{k-1-i}^{-1},\eta_{i+1}^{-1}\right\}\right)
= \rho(\xi_{i+1})
$$
for $i=0,\ldots,k-2$. The $L$-values of $\psi$ are the same as
those for $\wpsi$, but with the values at $i+1$ with $i\equiv
0,k-2$ (mod $p-1$) excluded. Moreover, $\xi_j = 0$ for even $j$,
so for $2\le i\le k-2$ we have proved
$$
L(\psi,i) = \begin{cases}
\rho(\xi_i) &\hbox{if $i$ is odd};\\
0 &\hbox{otherwise}.
\end{cases}
$$
This completes the proof of Theorem 1.2.

\section{Sharifi's Conjecture}

We remark that the Hecke eigenvalues of  the parabolic cohomology
class associated to the image of $\varphi$ in
$H^1_{par}(\Gamma_0,\F_p(\omega^{2-k}))$ for $T_2$ and $T_3$
correspond to those of the semi-cusp form in characteristic $p$ of
weight $2$ and type $\omega^{2-k}$
$$
s_{2,\omega^{2-k}}=\sum_{n\geq 1}\sum_{d|n}\omega^{2-k}(n/d)dq^n.
$$

 Let $f$ be a weight $k\geq 2$ cusp form whose Fourier expansion is
 given by $f(z)=\sum_{n\geq 1}a_ne^{2\pi inz}.$
Let  $L(f,s)=\sum_{n\geq 1}a_n n^{-s}$  be the complex
$L$-function of $f.$ Then, $f$ gives rise to a class $\phi_f\in
H^1_c(SL_2(\Z),V_{k-2}(\C))$ given by
$$\phi_f((x)-(y))=\int_{y}^{x} f(z)(zX+Y)^{k-2}dz.$$
It is well known that $L(\phi_f,\alpha)=
\frac{(-1)^{\alpha-1}(\alpha-1)!}{(2\pi i)^{\alpha}}L(f,\alpha)$
for all integers $\alpha$ with $1\leq \alpha\leq k-1.$

We will denote by $\psi_f$ the parabolic cohomology class
associated to $\phi_f$ in $H^1_{par}(SL_2(\Z),V_{K-2}(\C))$.

 Now let $p$ be an irregular prime and choose $k,$ $2\leq k<2p,$ such that
 $p\bigm|{B_k\over k}.$ Assume hypothesis $H_k$ holds for the irregular pair
 $(p,k)$.  Let $f$ be a normalized weight $k$ newform of level $1$ and let $\cO_f$ be the ring
of integers of the number field $K_f$ generated by the fourier
coefficients of $f$.  Let $\wp$ be a place above $p$ and assume
$f\equiv G_k\;(mod\;\wp)$, where $G_k$ is the Eisenstein Series of
level $1$ and weight $k$:
$G_k=-\frac{B_k}{2k}+\sum_{n=1}^{\infty}\sigma_{k-1}(n)q^n,$ with
$\sigma_{k-1}(n)=\sum_{d|n,d>0}d^{k-1}.$
\medskip

\noindent {\bf Conjecture.}
(Sharifi)\cite{SharifiSlides}\cite{SharifiL-values} {\it
 Under the above hypotheses, there exists $\Omega\in \C^\times$
such that $\frac{L(\psi_f,i)}{\Omega}\in \cO_f$, for all odd $i$
in the range $3\le i \le k-3$ and at least one of these numbers is
non-zero modulo $\wp$ and
$$
\frac{L(\psi_f,i)}{\Omega} \equiv \rho(\xi_i) \hbox{\rm \ \ (mod
$\wp$)}
$$
for all odd $i$ in the range  $3\leq i\leq k-3$. }

We remark that the computations of McCallum and Sharifi (see
\cite{MS1}-Theorem 5.1,\cite {MS2})imply that the space
$H^+_{k,eis,\{2,3\}}$ defined in the introduction is
one-dimensional for all irregular pairs $(p,k)$ with $p < 10,000$.
In fact, McCallum and Sharifi only compute the equivalent of the
eigenvalue of $T_2$, so in this range, we even have the stronger
statement
$$
\dim_{\F_p}\left(H^+_{k,eis,\{2\}}\right) = 1.
$$

\vfill\eject

\end{document}